\documentclass{article}

\usepackage[ansinew]{inputenc}
\usepackage{graphicx}
\usepackage{amsmath}
\usepackage{amsthm}
\usepackage{amssymb}

\theoremstyle{plain}
\newtheorem{thm}{Theorem}

\newtheorem{cor}[thm]{Corollary}
\newtheorem{prop}[thm]{Proposition}
\newtheorem{conj}{Conjecture}
\newtheorem{defin}{Definition}

\theoremstyle{remark}
\newtheorem{remark}{Remark}
\newtheorem{obs}[remark]{Observation}

\def\conv{\operatorname{conv}}

\title{Depth of segments and circles through points \\ enclosing many points: a note\thanks{%
Partially supported by CAM grant S-0505/DPI/0235-02.
Part of this work was done while the author was visiting the Mathematical Sciences Research Institute.}.
}

\author{Pedro A. Ramos \hskip10mm Raquel Viaña \\
    \small Departamento de Matem\'aticas\\
    \small Universidad de Alcal\'a \\
    \small Alcal\'a de Henares, Spain\\
    \small {\tt [pedro.ramos|raquel.viana]@uah.es}
}

\begin{document}

\maketitle

\begin{abstract}
Neumann-Lara and Urrutia showed in 1985 that in any set of $n$ points in the plane in general position
there is always a pair of points such that any circle through them contains at least $\tfrac{n-2}{60}$ points.
In a series of papers, this result was subsequently improved till $\tfrac{n}{4.7}$, which is currently the best
known lower bound. In this paper we propose a new approach to the problem that allows us, by using known
results about $j$-facets of sets of points in $\mathbb{R}^3$, to give a simple proof of a somehow stronger result:
there is always a pair of points such that any circle through them has, both inside and outside,
at least $\tfrac{n}{4.7}$ points.
\end{abstract}

\section{Introduction}

The problem that we address in this work was proposed by Neumann-Lara and Urrutia in \cite{nlu},
where the following result is shown: given a set $P$ of $n$ points in the plane in general position --
no three of them are collinear and no four of them are cocircular -- there is always a pair of points
$p,q\in P$ such that every circle through $p$ and $q$ contains at least $\left\lceil\tfrac{n-2}{60}\right\rceil$
other points of $P$.
In a series of papers \cite{hrw,bssu,h} this bound was slightly improved and, shortly afterwards,
Edelsbrunner et al. \cite{ehss}, by using techniques related to the complexity of higher order Voronoi diagrams,
showed a bound of $(\tfrac{1}{2}-\tfrac{1}{\sqrt{12}})n+O(1) \approx \tfrac{n}{4.7}$, which is the best
currently known lower bound for the problem. Regarding the upper bound, in \cite{hrw} Hayward et al. constructed
a set of $4m$ points such that for any two of them there are circles passing through them and containing less than
$m$ points.
Therefore, this example shows that $\lceil\tfrac{n}{4}\rceil-1$ is an upper bound for the problem.
In the same paper, the authors study the problem for sets of points in convex position, and give a bound of
$\lceil\tfrac{n}{3}\rceil-1$, which is also shown to be tight. Urrutia \cite{u} has conjectured that $\tfrac{n}{4}$ is,
up to perhaps an additive constant, the tight bound for the general problem.

In this note we give an alternative proof of the result by Edelsbrunner et al., transforming the problem from
circles in the plane to planes in the space. We introduce the concept of depth of a segment in a set of points
$P\subset\mathbb{R}^3$ and, by using known results about the number of $j$-facets, we show that there is always
a pair of points such that every circle through them has, both inside and outside, at least $\tfrac{n}{4.7}$ points.
Furthermore, we propose a new conjecture about the maximal number of segments with depth $k$ that a set of points
in convex position can have, which implies a stronger version of the original conjecture.

\section{Transforming the problem}

We use the well known transformation which maps the point $p=(p_x,p_y)\in\mathbb{R}^2$ to the point
$\hat{p}=(p_x,p_y,p_x^2+p_y^2)\in\mathbb{R}^3$ in the paraboloid $z=x^2+y^2$. Among the useful properties of
this transformation (see, for instance, \cite{e}) we will use the next one:

\begin{obs}
Given three non collinear points $p,q,r\in\mathbb{R}^2$, a point $s$ is inside the circle through them if and only
if point $\hat{s}$ is bellow the plane defined by $\hat{p},\hat{q},\hat{r}\in\mathbb{R}^3$.
\end{obs}
%\begin{demostracion}
%La demostración.
%\end{demostracion}

Therefore, the original problem is transformed into this one: given a set of $n$ points in the paraboloid
$z=x^2+y^2$, show that there exist a pair of points such that any plane passing through them leaves bellow at least
$\lceil \tfrac{n}{4}\rceil - 1$ points. This motivates the following definition:
%In the following, we assume that $P$ is a set of points in the space in general position, that is to say,
%no four points are on a common plane.

\begin{defin}
Given a set of points $P\subset\mathbb{R}^3$ and two points $p,q\in P$, the {\em depth} of segment $pq$ is defined
as the smallest integer $k$ such that any plane through $p$ and $q$ has on each side at least $k$ points of $P$.
\end{defin}

We observe that segments with depth zero are the edges of the convex hull and we are interested in showing that any
set of points has segments with ``high depth''.

We recall that, given points $p,q,r\in P$, the (oriented) triangle
$pqr$ is a $j$-facet of $P$ if it has exactly $j$ points on the positive side of its affine hull. Therefore, if
$pqr$ is a $j$-facet, its edges have depth at most~$j$. A subset $T\subset P$ is a $k$-set if it has $k$ points
and the sets $T$ and $P\smallsetminus T$ can be separated by a plane.
The number of $j$-facets of a set of points in $\mathbb{R}^d$
is related to the number of $(j\pm d)$-sets and obtaining tight bounds for these quantities, even for
$d=2$, is a famous open problem. The number of $(\leq j)$-facets is much better understood. In order to state the result,
we need some notation.

Let $e_j(P)$ be the number of $j$-facets of $P$ and let $E_j(P)=\sum_{i=0}^j e_i(P)$ be the number of
$(\leq j)$-facets. In \cite{w} Welzl shows the following:

\begin{thm}%[\cite{w}]
Let $P\subset\mathbb{R}^3$ be a set of $n$ points in general position. Then,
$$
E_j(P)\leq 2\Big[\binom{j+2}{2}\,n - 2\,\binom{j+3}{3} \Big] \qquad \text{if $\,\,0\leq 2j \leq n-4$.}
$$
Furthermore, the bound is tight and is achieved if and only if the set $P$ is in convex position.
\end{thm}

Because for a set of points in convex position $E_j(P)$ is known, the following result follows immediately:

\begin{cor} \label{cor:ej}
Let $P\subset\mathbb{R}^3$ be a set of $n$ points in convex position. Then,
$$
e_j(P)=E_j(P)-E_{j-1}(P)=2(j+1)n-2(j+1)(j+2) \qquad \text{if $\,\,0\leq 2j \leq n-4$.}
$$
\end{cor}

Next we use this result to bound the number of segments with depth at most $j$ for a set of points in
convex position. We denote by $s_j(P)$ the number of segments of $P$ with depth $j$ and by
$S_j(P)=\sum_{i=0}^j s_i(P)$ the number of segments with depth at most $j$.

\begin{prop}\label{p:T_j}
Let $P\subset\mathbb{R}^3$ be a set of $n$ points in convex position. Then,
$$
S_j(P) \leq  3(j+1)n-3(j+1)(j+2) \qquad \text{if $\,\,0\leq 2j \leq n-4$.}
$$
\end{prop}

\begin{proof}

Let $j$ be such that $0\leq 2j \leq n-4$. We claim that if $pq$ is a segment with depth at most $j$, then it is
contained in at least two $j$-facets of $P$. In order to prove the claim, consider first the case when the depth
is smaller than $j$ and let $\pi$ be an oriented plane passing through $p$ and $q$ and having less than $j$ points
in the positive side (denoted $\pi^+$ in Figure~\ref{fig1}). Because in the negative side of $\pi$ there are more
than $\lceil\tfrac{n}{2}\rceil$ points, if we rotate the plane around $pq$ in a direction we find, before having
rotated 180º, a point $r$ such that the plane $\pi_1$ passing through $p$, $q$ and $r$ leaves on the positive side
exactly $j$ points of $P$ and, therefore, $pqr$ (oriented conveniently) is a $j$-facet of $P$. In the same way, if
we rotate plane $\pi$ in the opposite direction, we find another point $s$ and, thus, another $j$-facet containing
segment $pq$. Finally, if the depth of $pq$ is $j$, we observe that the first point that we find when the plane
rotates must be in the negative side of the plane and thus it defines a $j$-facet.

\begin{figure}
\begin{center}
\includegraphics{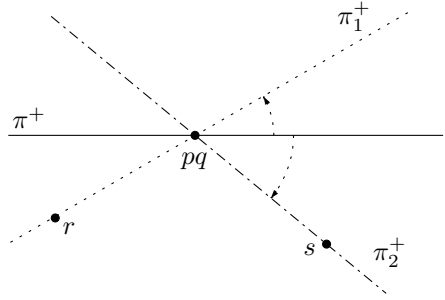}
\caption{Illustration for the proof of Proposition~\ref{p:T_j}.}
\label{fig1}
\end{center}
\end{figure}

Because each $j$-facet has 3 edges, it follows that $2S_j(P)\leq 3e_j(P)$ and, from Corollary~\ref{cor:ej} we get
$$
S_j(P) \leq \frac{3}{2}\, e_j(P) = 3(j+1)n-3(j+1)(j+2) \qquad \text{for $\,0\leq 2j \leq n-4$.}
$$
\end{proof}

We are ready to show the main result of this paper.

\begin{thm}
In a set $P\subset\mathbb{R}^3$ of $n$ points in convex position there exist segments with depth at least
$$
\Bigl(\frac{1}{2}-\frac{1}{\sqrt{12}}\Bigr)\,n + O(1)  \approx \frac{n}{4.7}.
$$
\end{thm}

\begin{proof}
Because $n$ determine $\binom{n}{2}$ segments, while $S_j(P)$ is smaller than $\binom{n}{2}$ there must be
segments with depth bigger than $j$. Therefore, from Proposition~\ref{p:T_j} we get
$$
3(j+1)n-3(j+1)(j+2) = \binom{n}{2},
$$
whose smaller solution is
$$
j=\frac{n-3}{2}-\Bigl( \frac{(n-2)^2-1}{12} \Bigr)^{1/2} = \Bigl(\frac{1}{2}-\frac{1}{\sqrt{12}}\Bigr)\,n + O(1).
$$

\end{proof}

Finally, if we apply this result to the original problem of circles passing through pairs of points, we
obtain immediately the following result:

\begin{cor}
Let $P$ be a set of $n$ points in the plane in general position. There always exists a pair of points $p,q\in P$ such
that every circle through $p$ and $q$ has, both inside and outside, at least
$$
\Bigl(\frac{1}{2}-\frac{1}{\sqrt{12}}\Bigr)\,n + O(1)  \approx \frac{n}{4.7}
$$
points of $P$.
\end{cor}

\section{A new conjecture}

We propose a new conjecture which has arisen during our study of this problem.

\begin{conj}\label{c2}
Let $P\subset\mathbb{R}^3$ be a set of $n$ points in convex position and let $s_j(P)$ be the number of
segments with depth $j$. Then,
$$
s_j(P) \leq 3n-8j-6 \qquad \text{if $\,\,0\leq j \leq \lceil\tfrac{n}{4}\rceil - 1$.}
$$
\end{conj}

Of course, the result is obvious (with equality) for $j=0$ and it is easy to give an almost
tight bound for $j=1$:

\begin{prop}
Let $P\subset\mathbb{R}^3$ be a set of $n$ points in convex position. Then,
$$
s_1(P) \leq 3n-12.
$$
\end{prop}

\begin{proof}
A segment $uv$ has depth one if and only if it is not an edge of the convex hull of $P$,
denoted by $\conv (P)$, but there exists a point
$p\in P$ such that $uv$ is an edge of $\conv(P\smallsetminus\{p\})$. If we denote by $\delta(p)$ the
number of vertices adjacent to $p$ in $\conv (P)$, the number of new edges in $\conv (P\smallsetminus\{p\})$
is exactly $\delta(p)-3$. Therefore,
\begin{equation}\label{eq:s1}
s_1(P) \leq \sum_{p\in P} (\delta(p)-3) = 3n-12.
\end{equation}
\end{proof}

\begin{remark}
The inequality in (\ref{eq:s1}) is strict if there is a segment $uv$ with depth one and points $p$ and $q$
such that $uv$ is an edge both of $\conv (P\smallsetminus\{p\})$ and $\conv (P\smallsetminus\{q\})$.
In this situation, we say that segment $uv$ is generated by two points.
It is easy to see that a segment with depth one cannot be generated by more than two points. Therefore,
Conjecture~\ref{c2} for $s_1(P)$ is equivalent to show that there are always at least two segments generated
by two points.
% (these segments are $r_2s_2$ and $q_2s_{m-1}$ in Figure~\ref{fig2}.b).
\end{remark}

\begin{figure}
\begin{center}
\includegraphics[width=0.8\textwidth]{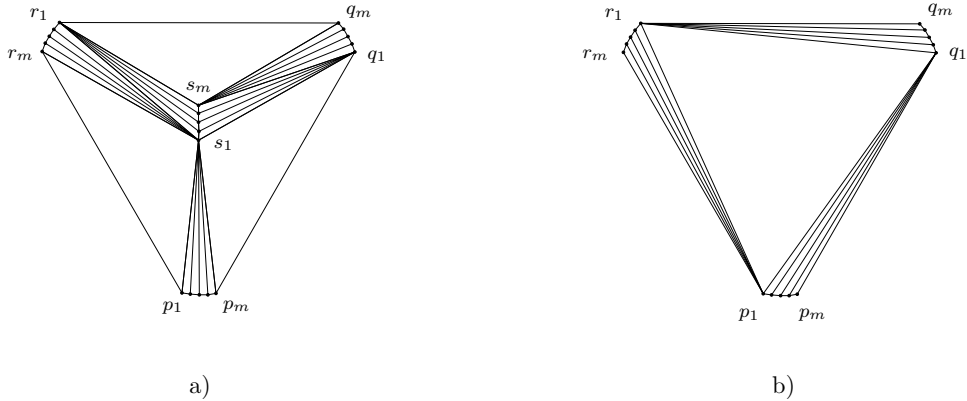}
\caption{Construction reaching $s_j(P) = 3n-8j-6$ for
$j=0,\ldots,\tfrac{n}{4} - 1$.}
\label{fig2}
\end{center}
\end{figure}

In the following we construct a set $P\subset\mathbb{R}^3$ such that $s_j(P) = 3n-8j-6$ for
every $j=0,\ldots,\tfrac{n}{4} - 1$, thus showing that the bound in Conjecture~\ref{c2}
would be tight. Consider the arc of circle $C=\{(x,y,z)\in\mathbb{R}^3\,|\,x^2+z^2=1,y=0,x>0.99\}$
and rotate it $45º$ counterclockwise around the $x$ axis. Let $n=4m$, put points $C_p=\{p_1,\ldots,p_m\}$
in $C$ and perturb them slightly to achieve general position. Now construct points $C_q$ and $C_r$
by rotating $C_p$ around the $z$ axis, $120º$ and $240º$, respectively. Finally, consider the arc
$C'=\{(x,y,z)\in\mathbb{R}^3\,|\,x^2+z^2=1,y=0,z>0.99\}$ and put the rest of the points,
$C_s=\{s_1,\ldots,s_m\}$, near $C'$ but slightly perturbed to achieve general position. The convex
hull of $P=C_p\cup C_q \cup C_r\cup C_s$ is shown in Figure~\ref{fig2}.a (top view) and Figure~\ref{fig2}.b
(bottom view).

The fact that $s_j(P)=4n-8j-6$ for $j=0,\ldots,\tfrac{n}{4}-1$ can be easily checked taking into
account the following simple observations:
\begin{itemize}
\item[--] A segment $s$ has depth $j$ if it is in the convex hull of $P\smallsetminus T$ for
some $j$-set $T$ and it is not in the convex hull of $P\smallsetminus S$ for any $k$-set $S$
with $k<j$.
\item[--] Given $T\subset P$ with $|T|<n/4$, the convex hull of $P'=P\smallsetminus T$ has
``the same structure'' as $\conv (P)$, i.e., consecutive points in each of the chains are adjacent,
the first point in $C_s'$ is adjacent to all the points in $C_r'$ and $C_p'$, and so on.
\end{itemize}

We conclude the note stating a direct implication of the previous conjecture.
Because
$$
\sum_{j=0}^{\lfloor\tfrac{n}{4}\rfloor -2} (3n-8j-6) \leq \binom{n}{2}-(n+2),
$$
Conjecture~\ref{c2} would imply:
\begin{conj}\label{c3}
For every set of $n$ points in the plane in general position, there are always $n+2$ pairs
of points such that any circle through them has, both inside and outside,
at least $\lfloor\tfrac{n}{4}\rfloor - 1$ points.
\end{conj}

%As an example, consider a set of points in the moment curve.
%
%\begin{prop}
%Let $t_1,t_2,\ldots,t_n$ be distinct real numbers, $p_i=(t_i,t_i^2,t_i^3)$ and $S=\{p_1,p_2,\ldots,p_n\}$ a set of
%points in the moment curve. Then
%$$
%\delta_j(S) \leq 3n-9j-6 \qquad \text{if $0\leq j \leq \tfrac{n}{4}$.}
%$$
%\end{prop}
%
%\begin{proof}
%
%\end{proof}
%
%Obsérvese que esta última conjetura implicaría que hay ``muchos'' pares de puntos como el propuesto por la
%conjetura original. En concreto, para el caso en que $n$ es múltiplo de 4, se tiene que
%$$
%\sum_{j=0}^{\frac{n}{4}-2} (3n-8j-6) = \binom{n}{2} - (n+2),
%$$
%de donde se deduciría que hay $n+2$ pares de puntos tales que todas las circunferencias que pasan por ellos tienen,
%tanto dentro como fuera, al menos $\tfrac{n}{4}-1$ puntos.
%
%Por último, merece la pena mencionar que el mismo ejemplo de \cite{hrw} (levantado al paraboloide) que muestra que
%$\lceil\tfrac{n}{4}\rceil-1$ es una cota superior para el problema original demuestra también que la cota para el
%número de segmentos de profundidad $j$ propuesta en la Conjetura~\ref{c2} es ajustada para
%$j=0,\ldots,\lceil\tfrac{n}{4}\rceil-1$.

\section{Acknowledgements}

I would like to thank Julian Pfeiffle for his constructions using Polymake and Boris Aronov, Imre Bárány, David Orden,
and Micha Sharir for helpful discussions.

\end{document}